\input amstex
\documentstyle {amsppt}
\UseAMSsymbols \vsize 18cm \widestnumber\key{ZZZZZ}

\catcode`\@=11
\def\displaylinesno #1{\displ@y\halign{
\hbox to\displaywidth{$\@lign\hfil\displaystyle##\hfil$}&
\llap{$##$}\crcr#1\crcr}}
\def\ldisplaylinesno #1{\displ@y\halign{
\hbox to\displaywidth{$\@lign\hfil\displaystyle##\hfil$}&
\kern-\displaywidth\rlap{$##$} \tabskip\displaywidth\crcr#1\crcr}}
\catcode`\@=12

\refstyle{A}

\let \ol=\overline

\let \ti=\widetilde

\font\sc=cmcsc8 at 10pt

 \font\srm=cmr10 at 7.5pt

\font\main=cmsy10 at 10pt

\font \fin=lasy8 at 15.4 pt

\def \A{\mathop{\hbox{\main A}}\nolimits}

\topmatter
\title The value of the global intertwining operators on spherical vectors\endtitle

\rightheadtext{}
\author Volker Heiermann\endauthor
\address Aix-Marseille Universit\'e, CNRS, Centrale Marseille, I2M, UMR 7373, 13453 Marseille, France;
\endaddress
\email volker.heiermann\@univ-amu.fr \endemail

\thanks The author has benefitted from a grant of Agence
Nationale de la Recherche with reference  ANR-13-BS01-0012 FERPLAY.
\endthanks

\abstract
Let $F$ be a global field, $G$ an unramified quasi-split reductive group over $F$ and $\chi $ an everywhere unramified automorphic character of a maximal maximally split torus of $G$. Using Langlands-Shahidi theory, we compute the meromorphic function defined by the action of a global standard intertwining operator associated to $\chi $ on a spherical vector and show that the ratio of its poles in the positive Weyl chamber is well behaved.
\endabstract

\endtopmatter
\document

Let $F$ be a global field and let $G$ be the group of $F$-rational points of a quasi-split connected reductive group that splits over an extension $E$ of $F$, which is unramified in each non-Archimedean place (in the number field case, we do not assume that the Archimedean places split). If $v$ is a place of $F$, denote by $G_v$ the group of its $F_v$-rational points.

Fix a Borel subgroup $B$ of $G$ defined over $F$ and denote by $U$ its unipotent radical and by $T$ a Levi subgroup that is a maximally split, but possibly non split, torus. Denote by $T_0$ the maximal split subtorus of $T$, by $w$ an element in its Weyl group and by $\ti{w}$ a representative of $w$ in $G$ chosen as in \cite{Sh90}. (Inside the proofs, $w$ will later also denote occasionally a place of an extension of $F$ - we hope that this will not confuse the reader.)

We denote by $\Sigma _{red}$ the set of reduced roots of $T_0$ in the Lie algebra of $G$ and by $\Sigma_{red} (B)$ the subset of the roots which are positive w.r.t. $B$. The elements of $\Sigma_{red} (B)$ will sometimes be called relative positive roots.

Fix at all non-Archimedean places $v$ a hyperspecial compact subgroup $K_v$ and at the Archimedean places a special compact subgroup $K_v$, such that $G_v=B_vK_v$ at all places. Put $\Bbb K= \prod _v K_v$. Let $\Bbb A=\Bbb A_F$ be the ring of ad\`eles associated to $F$. Then, $G(\Bbb A)$ is the restricted tensor product of the $G_v$ w.r.t. the $K_v$, and the equality $G(\Bbb A)=B(\Bbb A)\Bbb K$ holds. We assume in addition that choices have been made so that for each $F$-parabolic subgroup $P=MV$ containing $B$, one has $P(\Bbb A)\cap\Bbb K=(M(\Bbb A)\cap\Bbb K)(V(\Bbb A)\cap\Bbb K)$ with $M(\Bbb A)\cap\Bbb K$ maximal compact. One gives $\Bbb K$ the measure $1$, and, if $V$ is a unipotent group over $F$, then the Haar measure on $V(\Bbb A)$ is normalized such that the quotient $V(F)\backslash V(\Bbb A)$ gets measure $1$, where one takes on $V(F)$  the counting measure. If one fixes a measure on $T(\Bbb A)$, then this defines a measure on $G(\Bbb A)$ \cite{MW, I.1.13}.

Fix $\chi=\otimes_v\chi _v$ a unitary everywhere unramified automorphic character of $T(\Bbb A)$. (This means that $\chi _v$ factors also at the Archimedean places  through the absolute value.)

If $\alpha $ is a positive relative root for $G$, we will denote by $M_{\alpha }$ the Levi subgroup of $G$ of relative semi-simple rank $1$ associated to $\alpha $. Denote by $E_{\alpha }$ the minimal splitting field for the derived group $M_{\alpha }^{der}$ and by $T_{\alpha }$ the maximally split torus of $M_{\alpha }^{der}$ contained in $T$. We will denote by $\underline{\alpha}^\vee $ the $F$-rational map $Res_{E_{\alpha }/F}\Bbb G_m\rightarrow T_{\alpha }$ defined in the following way: let $\Gamma _F$ be the absolute Galois group of $F$ with respect to a given separable closure $F^{sep}$ of $F$, $\Gamma _{E_{\alpha }}$ the absolute Galois group of $E_{\alpha }$ seen as a subgroup of $\Gamma _F$ and $\Gamma _{E_{\alpha }\backslash F}$ a set of representatives of $\Gamma _{E_{\alpha }}\setminus\Gamma _F$. Fix an absolute simple root $\underline{\alpha}_1$ for $M_{\alpha }$ that restricts to $\alpha $ and denote by $\underline{\alpha}_1^{\vee }$ the corresponding (absolute) coroot. Then, we put $\underline{\alpha}^\vee =((\underline{\alpha}_1^{\vee })^{\sigma })_{\sigma\in\Gamma _{E_{\alpha }\backslash F}}$. More precisely, over $F^{sep}$, one has $\underline{\alpha}^\vee((x_{\sigma })_{\sigma\in\Gamma _{E_{\alpha }\backslash F}})=\prod _{\sigma\in\Gamma _{E_{\alpha }\backslash F}}(\underline{\alpha}_1^{\vee })^{\sigma }(x_{\sigma })$ from which follows that $\underline{\alpha}^\vee$ is $\Gamma_F$-invariant, i.e. defined over $F$. If $v$ is a place of $F$, then $\underline{\alpha}^\vee_v$ will denote $\underline{\alpha}^\vee$ over the $F_v$-rational points, i.e. $\underline{\alpha}^\vee_v: (Res_{E_{\alpha }/F}\Bbb G_m)(F_v)\rightarrow T_{\alpha }(F_v)$. In particular, $\chi\circ\underline{\alpha}^\vee$ will stand for the automorphic character $\otimes_v(\chi _v\circ\underline{\alpha}^\vee_v)$.

Denote by $X^*(T)$ the group of $F$-rational characters of $T$ and put $a_{T,\Bbb C}^*=X^*(T)\otimes\Bbb C$. One defines a map from $a_{T,\Bbb C}^*$ into the set of unramified automorphic characters of $T(\Bbb A)$ by sending $\alpha\otimes s$ to the character $t=(t_v)_v\mapsto\prod _v\vert\alpha(t_v)\vert_v^s$. It is a group homomorphism that is injective if $F$ is a number field and with kernel ${2\pi i\over\log q}\Bbb Z$ with $q$ equal to the cardinality of the ground field if $F$ is a function field. If $\lambda\in a_{T,\Bbb C}^*$, we will denote by $\chi _{\lambda }$ the product of the character $\chi $ by the unramified automorphic character of $T(\Bbb A)$ corresponding to $\lambda $.

Let $i_{B(\Bbb A)}^{G(\Bbb A)}$ be the functor of normalized parabolic induction. For $\lambda\in a_{T,\Bbb C}^*$, denote by $M(\lambda, \chi, \ti{w})$ the global standard intertwining operator $i_{B(\Bbb A)}^{G(\Bbb A)}\chi_{\lambda }\rightarrow i_{B(\Bbb A)}^{G(\Bbb A)}w\chi_{\lambda }$ defined, if $\lambda $ lies in the positive Weyl chamber corresponding to $B$, by
$$(M(\lambda, \chi, \ti{w})f)(g)=\int _{U(\Bbb A)\cap w\ol{U}(\Bbb A)w^{-1}}f( \ti{w}^{-1}ug)du,\eqno{(1)}$$
and by $\langle\cdot,\cdot\rangle$ the $G(\Bbb A)$-invariant standard natural pairing between $i_{B(\Bbb A)}^{G(\Bbb A)}\chi_{\lambda }$ and $i_{B(\Bbb A)}^{G(\Bbb A)}\chi_{-\ol{\lambda }}$ defined by the measure on $G(\Bbb A)$ (cf. \cite{MW, II.1.6}). It is semi-linear in the first variable and linear in the second variable. If $H$ is a linear algebraic group over $F$, $H^{der}$ will denote its derived group and $\ti{H}$ its simply connected covering group over $F$, so that $\ti{H^{der}}$ is the simply connected covering of the derived group of $H$ over $F$.

The main result of this note is the following theorem, which until now has only been stated in special cases.

\null {\bf Theorem:} \it  Let $\chi=\otimes_v\chi _v$ be a unitary everywhere unramified automorphic character of $T(\Bbb A)$.
Let $\phi _1$ (resp. $\phi _1^{\vee }$) be the $\Bbb K$-spherical vector in $i_{B(\Bbb A)}^{G(\Bbb A)}\chi_{\lambda }$ (resp. $i_{B(\Bbb A)}^{G(\Bbb A)}w\chi_{-\ol{\lambda }}$) with $\phi _1(1)=1$ (resp. $\phi _1^{\vee }(1)=1$).

One has $$\langle \phi _1^\vee , M(\lambda , \chi, \ti{w})\phi _1\rangle = \prod_{\alpha\in\Sigma_{red}(B)\cap w^{-1}\Sigma_{red}(\ol{B})}r_{\alpha }(\langle\lambda,\alpha^{\vee }\rangle,\chi\circ {\underline{\alpha}}^{\vee }),$$
where, for $s\in\Bbb C$ and $\eta $ an automorphic character,
$$\eqalign{&r_{\alpha }(s,\eta )\cr =&\cases{L_{F_{\alpha }}({s\over d_{\alpha }},\eta )\over \epsilon_{F_{\alpha }}({s\over d_{\alpha }},\eta )L_{F_{\alpha }}(1+{s\over d_{\alpha }}, \eta )},& if \ti{M_{\alpha}^{der}}\simeq Res_{F_{\alpha }/F}SL_2; \cr {L_{E_{\alpha }}({s\over 4d_{\alpha }},\eta )\over \epsilon_{E_{\alpha }}({s\over 4d_{\alpha }} ,\eta )L_{E_{\alpha }}(1+{s\over 4d_{\alpha }} , \eta )}\times\cr \qquad\times{L_{F_{\alpha }}({s\over 2d_{\alpha }} ,\eta\eta_{E_{\alpha }/F_{\alpha }})\over \epsilon_{F_{\alpha }}({s\over 2d_{\alpha }} ,\eta\eta_{E_{\alpha }/F_{\alpha }} )L_{F_{\alpha }}(1+{s\over 2d_{\alpha }} ,\eta\eta_{E_{\alpha }/F_{\alpha }} )},& if \ti{M_{\alpha}^{der}}\simeq Res_{F_{\alpha }/F}SU(2,1).\cr\endcases}$$
Here, $E_{\alpha }$ denotes the unramified quadratic extension of $F_{\alpha }$ over which $SU(2,1)$ splits, $\eta_{E_{\alpha }/F_{\alpha }}$ the quadratic character that defines $E_{\alpha }/F_{\alpha }$, $d_{\alpha }=[F_{\alpha }:F]$ and $\eta $ is an automorphic character of $(Res_{F_{\alpha }/F}\Bbb G_m)(\Bbb A_F)$ in the first case and of $(Res_{E_{\alpha }/F}\Bbb G_m)(\Bbb A_F)$ in the second one, while $L_{F_\alpha }$, $L_{{E_\alpha }}$ (resp. $\epsilon_{F_{\alpha }}$, $\epsilon_{E_{\alpha }}$) denote the (completed) Hecke $L$-functions (resp. $\epsilon $-factors) w.r.t. the fields $F_{\alpha }$, $E_\alpha$.

\null\it Remarks: \rm (i) If $\ti{M_{\alpha}^{der}}\simeq Res_{F_{\alpha }/F}SL_2$, the splitting field of $M_{\alpha }$ is $F_{\alpha }$ and, for $v$ a place of $F$, $(Res_{F_{\alpha }/F}SL_2)(F_v)=\prod_{w\vert v}SL_2(F_{\alpha ,w})$. Then, $\underline{\alpha}^{\vee }_v$ is the product of the usual coroots $F_{\alpha ,w}^*\rightarrow T(F_v)$ with $w$ going over the places of $F_{\alpha }$ over $v$.

If $\ti{M_{\alpha}^{der}}\simeq Res_{F_{\alpha }/F}SU(2,1)$, then the splitting field of $M_{\alpha }$ is $E_{\alpha }$, and one has
$Res_{F_{\alpha }/F}SU(2,1)(F_v)=\prod _{v'\vert v}SU(2,1)(F_{\alpha ,v'})$. Then, $\underline{\alpha}^{\vee }_v$ is a product of coroots $E_{\alpha ,w}^*\rightarrow T(F_v)$ with $w$ going over the places of $E_{\alpha }$ over $v$. The nature of the coroot differs if $w$ is invariant by the Galois group of $E_{\alpha }/F_{\alpha }$ or not. Denote by $v'$ the place of $F_{\alpha }$ below $w$. If $w$ is Galois invariant, then $SU(2,1)$ does not split over $F_{\alpha ,v'}$ and the coroot is given by the relative reduced positive coroot. Otherwise, $SU(2,1)(F_{\alpha ,v'})\simeq SL_3(E_{\alpha ,w})$ and there are two places $w_1$ and $w_2$ over $v'$. The coroots on $E_{\alpha ,w_1}^*$ and $E_{\alpha ,w_2}^*$ are given respectively by the two positive simple coroots for $SL_3$.

(ii) In $L_{F_{\alpha }}(s,\eta )$, it is the restriction $\eta _{\vert \Bbb A^{\times }_{F_{\alpha }}}$ that is considered. If $\eta =(\eta _w)_w$ with $w$ going over the places of $E_{\alpha }$, then $\eta _{\vert \Bbb A^{\times }_{F_{\alpha }}}=(\eta _{v'})_{v'}$ with $\eta _{v'}={\eta_w}_{\vert F_{\alpha, v'}}$ if $v'$ does not split in $E_{\alpha }$ and $\eta _{v'}=\eta_{w_1}\eta_{w_2}$ if $v'$ splits and $w_1$ and $w_2$ are the two places dividing $v'$, recalling that the restriction to $F_{\alpha ,v}$ is the diagonal embedding into $E_{\alpha , w_1}\oplus E_{\alpha , w_2}$ and that $E_{\alpha ,w_1}\simeq E_{\alpha ,w_2}\simeq F_{\alpha ,v}$ over $F$.

(iii) During the proof, we will use the Langlands-Shahidi normalization of the local standard intertwining operators. We want to stress that these are expressed in terms of the local automorphic $L$- and $\epsilon $-factors attached to the dual $\ti{r}$ of the adjoint representation $r$ of $^LT$ on $Lie(\ ^LU)$. This is imposed in the non-Archimedean case by Langlands' notation and convention, which have been taken up by Shahidi in his work. On the other hand, as it was pointed out to us by Shahidi, if one adopts the convention in \cite{Bo79}, one no longer needs the contragredients mentioned above as it is reflected in the formula for the constant term M(s) in page 53 of \cite{Bo79}. This is the case for the formulas for the local coefficients in the Archimedean setting in [Sh85], but to be consistent, the Langlands' convention has also to be used in the Archimedean case taken up by Shahidi in his work (see also \cite{Sh07, p. 310}).
\null

\it Proof: \rm By the multiplicity formula for global intertwining operators, it is enough to prove this for $G$ of relative semisimple rank $1$ over $F$. We have then a single reduced root $\alpha $ and can replace $\lambda $ by $s\ti{\alpha }$, where $\ti{\alpha }$ is as defined in \cite{Sh90, p. 278, l.-3}, so that $\langle\lambda,\alpha^{\vee }\rangle=d_{\alpha }s$ if $\ti{G^{der}}\simeq Res_{F_{\alpha }/F}SL_2$ and $\langle\lambda,\alpha^{\vee }\rangle=4d_{\alpha }s$ if $\ti{G^{der}}\simeq Res_{F_{\alpha }/F}SU(2,1)$ \cite{Sh10, p.10}, and the equalities above can be expressed in terms of the complex variable $s$. Denote by $w_{\alpha }$ the simple reflection associated to $\alpha $ and write in the sequel $F'$ for $F_{\alpha }$ and $E'$ for $E_{\alpha }$. (Remark that $E'=E$, as $G$ is of relative rank $1$.) We will \it not \rm (and do not need to) assume here that $w_{\alpha }$ fixes $\chi $.

Fix $\psi _F$ a non trivial unitary character of $F\backslash\Bbb A_F$, $\psi _F=\otimes_v\psi _{F_v}$.

For a place $v$ of $F$, denote by $A_v(s\ti{\alpha }, \chi_v, \ti{w}_{\alpha })$ the corresponding local standard intertwining operator defined for $\Re(s)>0$ by the converging integral analog to $(1)$ and w.r.t. the measure on $U_v$ derived from the measure on $F_v$, which is self-dual w.r.t. $\psi _{F_v}$. (This is in accordance to Langlands-Shahidi theory \cite{Sh10, p. 133}, \cite{Lo15}.) The measure on $U(F)\backslash U(\Bbb A_F)$ has then the correct properties and, in particular, one has $M(s\ti{\alpha }, \chi, \ti{w}_{\alpha })=\otimes_vA_v(s\ti{\alpha }, \chi_v, \ti{w}_{\alpha })$ \cite{Sh10, 6.3.3}.

Suppose first $\ti{G^{der}}$ of relative rank $1$ over $F_v$. Then, $\alpha $ is the reduced $F_v$-relative positive root for $T$ in $B$. Denote by $r_{\alpha }$ the adjoint representation of $^LT$ on the Lie algebra of $^LU$, by $r_{\alpha ,j}$ its irreducible components ordered as in \cite{Sh90}, and by $r_{\alpha,j}^{\vee }$ the contragredient representation. The normalized intertwining operator in $v$ is then given in the number field case by (cf. \cite{Sh90, (7.12)}\footnotemark\footnotetext{The expression given here for the normalizing factor follows from \cite{Sh90, (7.12)} by applying the compatibility with unramified twist, which is a consequence of identity \cite{Sh90, (3.12)}  for the $\gamma $-factor, remarking that $s+s_0$ has to be replaced therein by $s-is_0$ as confirmed by the author. \hfill
In fact, the following paragraph was prepared by the author of [Sh90] in response to inquiries by the author concerning the sign of this shift: To explain the shift $s- is_0$ in $\gamma_i$, we note that while multiplication by $i$ is rather obvious, the shift in sign is more delicate as we now detail: The negative sign is a reflection of the appearance of the contragredient of the representation as opposed to representation itself in the right hand side of equation (3.11) in [Sh90]. Inductively we may assume $m=1$. We have
$$C_\chi(s,\widetilde\pi(s_0))=C_\chi(s+s_0,\widetilde\pi)=\gamma(s+s_0,\pi)=\gamma(s-(-s_0),\pi).\eqno{(0)}$$
On the other hand
$$C_\chi(s,\widetilde\pi(s_0))=\gamma(s,\widetilde{\widetilde\pi(s_0)})=\gamma(s,\pi(-s_0)),$$
which by (0) equals:
$$ =\gamma(s-(-s_0),\pi).$$
Therefore
$$\gamma(s,\pi(-s_0))=\gamma(s-(-s_0),\pi).$$
Now, changing $s_0$ to $-s_0$, we have
$$\gamma(s,\pi(s_0))=\gamma(s-s_0, \pi).$$}, \cite{Sh10, (8.5.1)})
$$\A _v(s\ti{\alpha }, \chi_v, \ti{w}_{\alpha })=\prod_{j=1}^m {\epsilon _v(js,\chi _v,r^{\vee }_{\alpha ,j},\psi _{F_v})L_v(1+js,\chi _v,r^{\vee }_{\alpha ,j})\over L_v(js,\chi _v,r^{\vee }_{\alpha ,j})} A_v(s\ti{\alpha }, \chi_v, \ti{w}_{\alpha }).$$ This makes also sense in the function field case, as we only consider unramified characters and the local factors are given by local class field theory \cite{Sh90, 3.5.1}.
In addition, the local factors are equal to the ones relative to the simply connected covering of the derived group of $G$ as the adjoint action is trivial on the center of the connected component of $^LG$ and, consequently, the representation $r_{\alpha }$ factors through the image of $^LT$ in the L-group of $\ti{G^{der}}$.

Remark that $Res_{F'/F}$ is, as functor of points on algebraic groups over $F_v$,  isomorphic  to $\prod_{v'\vert v}Res_{F'_{v'}/F_v}$, the product being taken over the places $v'$ of $F'$ that divide $v$. The group $Res_{F'_{v'}/F_v}SU(2,1)$ is of relative rank $1$, if and only if $SU(2,1)$ does not split over $F'_{v'}$ (or equivalently $v'$ does not split in $E'$). Otherwise, $Res_{F'_{v'}/F_v}SU(2,1)$ is isomorphic to $Res_{F'_{v'}/F_v}SL_3$.

When $\ti{G^{der}}\simeq_{F_v}Res_{F'_{v'}/F_v}SL_2$ with $v'$ a place of $F'$ over $v$, then $m=1$, the local $L$-factor is the Hecke one relative to $\chi_v\circ\underline{\alpha}^{\vee }_v$ and the field $F'_{v'}$ and the local $\epsilon $-factor the Hecke one relative to the additive character $\psi _{F_v}\circ Tr_{F'_{v'}/F_v}$ multiplied by $\lambda(F'_{v'}/F_v,\psi _{F_v})$ (cf. \cite{KeSh88, (2.9) and (2.10)} for the non-Archimedean case, that agrees with \cite{Lo, 1.4} in the function field case, and \cite{Sh85, p. 989-990} \footnotemark\footnotetext{To be consistent and get $\ti{r}$ also in the Archimedean case (compare also with remark (iii) above), one has to use equation \cite{KeSh88, (2.1)} in \cite{Sh85} also in the Archimedean case. Remark that in the original context of \cite{Sh85, theorem 3.1}, $C(-2s\rho_\theta,\dots)$ has to be replaced by $C(2s\rho_\theta,\dots)$.} for the Archimedean case). When $\ti{G^{der}}\simeq_{F_v} Res_{F'_{v'}/F_v}SU(2,1)$ and $SU(2,1)$ does not split over $F'_{v'}$, then $m=2$, the local $L$ factors for $r_{\alpha ,1}^{\vee }$ are those of Hecke relative to $\chi_v\circ\underline{\alpha}^{\vee }_v$ and the composite field $E'F'_{v'}$, the local $L$-factors for $r_{\alpha ,2}^{\vee }$ are those of Hecke relative to $\eta_{E'F_v/F_v}({\chi_v}_{\vert F_v^*})\circ\underline{\alpha}^{\vee }_v$ and the field $F'_{v'}$. The $\epsilon $-factors are the appropriate ones multiplied respectively by $\lambda(F'_{v'}/F_v,\psi _{F_v})$ and $\lambda(E'F'_{v'}/F_v,\psi _{F_v})$ (cf. \cite{KeSh88, (2.13) and (2.14)} for the non-Archimedean case, that agrees with \cite{Lo, 1.4} in the function field case, and \cite{Sh85, p. 991-992}$^{\srm 2}$ for the Archimedean case).

If $\ti{G^{der}}\simeq_{F_v} Res_{F'_{v'}/F_v}SU(2,1)$ and $SU(2,1)$ splits over $F'_{v'}$, then $\ti{G^{der}}\simeq_{F_v} Res_{F'_{v'}/F_v}SL_3$, $\ti{w}_{\alpha }$ becomes a representative of the longest element of the Weyl group of $Res_{F'_{v'}/F_v}SL_3$. Denote the two simple positive relative coroots by $\alpha _1^{\vee }$ and $\alpha _2^{\vee }$. From the product formula for the normalized intertwining operator, it follows that
$$\eqalign{&\A _v(s\ti{\alpha }, \chi_v, \ti{w}_{\alpha })\cr =&{\epsilon _{v'}(s,\chi_v\circ\alpha_1^{\vee },\psi _{F'_{v'}})L_{v'}(1+s,\chi_v\circ\alpha_1^{\vee })\over L_{v'}(s,\chi_v\circ\alpha_1^{\vee })}{\epsilon _{v'}(s,\chi_v\circ\alpha_2^{\vee },\psi _{F'_{v'}})L_{v'}(1+s,\chi_v\circ\alpha_2^{\vee })\over L_{v'}(s,\chi_v\circ\alpha_2^{\vee })}\times\cr &\times{\epsilon _{v'}(2s,\chi _v\circ(\alpha_1^{\vee }+\alpha _2^{\vee }),\psi _{F'_{v'}})L_{v'}(1+2s,\chi _v\circ(\alpha_1^{\vee }+\alpha _2^{\vee }))\over L_{v'}(2s,\chi _v\circ(\alpha_1^{\vee }+\alpha _2^{\vee }))} A_v(s\ti{\alpha }, \chi_v, \ti{w}_{\alpha }),\cr}$$ where the local factors are once more the Hecke one's relative to $F'_{v'}$, the $\epsilon $-factors being taken relative to the additive character $\psi _{F_v}\circ Tr_{F'_{v'}/F_v}$ and multiplied by $\lambda(F'_{v'}/F_v,\psi _F)$ \cite{KeSh88, 2.10}.

In general, if the simply connected covering of the derived group of $G$ is isomorphic over $F_v$ to a product of semi-simple simply connected groups of relative rank $1$ or of the type $Res_{F'_{v'}/F_v}SL_3$, then the normalized intertwining operator is the tensor product of the ones explained above.

Writing $\A(s\ti{\alpha }, \chi, \ti{w}_{\alpha })=\otimes_v\A_v(s\ti{\alpha }, \chi_v, \ti{w}_{\alpha })$ and $r_v(s,\chi_v\circ\underline{\alpha}^{\vee }_v,\psi _{F_v})$ for the reciprocal of the normalizing factor so that
$$\A _v(s\ti{\alpha }, \chi_v, \ti{w}_{\alpha }) =r_v(s,\chi_v\circ\underline{\alpha}^{\vee }_v,\psi _{F_v})^{-1}A_v(s\ti{\alpha }, \chi_v, \ti{w}_{\alpha }),$$
everything comes down to show that $r_v(s,\chi_v\circ\underline{\alpha}^{\vee }_v,\psi _{F_v})$ is the product, over the places $v'$ of $F'$ over $v$, of the Euler-factors of $r_{\alpha }(s,\chi\circ\alpha^{\vee })$ defined by the choice of the character $\psi _F\circ Tr_{F'/F}$, and that $\langle \phi_1^{\vee }, \A(s\ti{\alpha }, \chi, \ti{w}_{\alpha }) \phi_1\rangle =1.$

For the first point, the Euler product expression, as $Res_{F'/F}$ is over $F_v$ isomorphic to $\prod_{v'\vert v}Res_{F'_{v'}/F_v}$, things are immediate when $\ti{G^{der}}\simeq Res_{F'/F}SL_2$ or when $\ti{G^{der}}\simeq Res_{F'/F}SU(2,1)$ and $SU(2,1)$ does not split over any of the $F_{v'}'$. When $\ti{G^{der}}\simeq Res_{F'/F}SU(2,1)$ and some of the places $v'$ over $v$ split in $E'$, then one has to check that the Euler factor corresponding to such a $v'$ equals the reciprocal of the normalizing factor relative to $SL_3(F_{v'})$, remarking that this factor does not depend on the isomorphism class of $F'_{v' }$. But, this follows from the remarks (i) and (ii) after the statement of the theorem and the above computation of the normalizing factor for $SL_3$.

Concerning the second point, the equality $\langle \phi _1^{\vee }, \A(s\ti{\alpha }, \chi, \ti{w_{\alpha }})\phi _1\rangle =1$, it will be enough to show that
$$(\A _v(s\ti{\alpha }, \chi_v, \ti{w_{\alpha }})\phi _{1,v})(1)=\lambda_{\psi_{F_v}}(\ti{w_{\alpha }}),\eqno{(2)}$$
where $\lambda_{\psi_{F_v}}(\ti{w_{\alpha }})$ is the constant defined in the non-Archimedean case in \cite{KeSh88, p. 80, (4.1)} (this agrees with \cite{Lo15, 2.5} in the function field case)  deduced from the  Langlands $\lambda$-function and by the same formula in the Archimedean case \cite{Sh10, p. 136}. In fact, the product over all places is one: this is remarked in \cite{KeSh88, p. 83, l.1} for the number field case. In fact, it is a direct consequence of the compatibility of the global $\epsilon $-factor with induction, which is also true in the function field case \cite{De73, 3.12 C}, remarking that the Langlands local factors are the Deligne ones for the self-dual additive measure \cite{Ta79 (3.6.4)}, see also \cite{De73, 5.6.2}.

We will distinguish the cases of $v$ non-archimedean, $v$ real archimedean and $v$ complex archimedean.

Suppose first that $v$ is a non-Archimedean place and that $F_v$ is of characteristic $0$. Then $G_v$ is a $p$-adic group. One knows that the action of the normalized intertwining operator $\A _v(s\ti{\alpha }, \chi_v, \ti{w_0})$ on the $K_v$-spherical vector $\phi_{1,v}$ does not depend on $s$ (\cite{Sh90, remark 7.10} with \cite{Ar89, Theorem 2.1}). So, to compute the value in (2), we can assume $s=0$ and $\chi_v=1$, as $\chi _v$ is unramified by assumption and $r_{\alpha ,v}(s,\chi _{s_0\ti{\alpha }},\psi _{F_v})=r_{\alpha ,v}(s+s_0,1,\psi _{F_v})$ \cite{Sh90, (3.12), (7.1), (7.2), (7.4)}. As the Whittaker functional is non-trivial on the spherical vector (cf. \cite{KeSh88}, end of proof of theorem {\bf 4.1})\footnotemark\footnotetext{In \cite{KeSh88}, it is assumed that $K_v=G(O_v)$, but, after an isomorphism, this holds for hyperspecial groups and even if $G$ is not defined on $O_v$, as hyperspecial maximal compact subgroups form a single orbit under the action of the adjoint group (cf. \cite{Tit, 2.5}) and an isomorphism over $F$ sends hyperspecial maximal compact subgroups to hyperspecial ones.}, one can apply \cite{KeSh88, 3.4} with the formula for $\lambda_{\psi_{F_v}}(\ti{w_{\alpha }})$ in \cite{KeSh88, p. 80, (4.1)} and the fact that $L(s,1,r_{\alpha })=L(s,1,r_{\alpha }^{\vee })$ by the identity (7.2) of \cite{Sh90}.

If $F_v$ is of positive characteristic, then the computation of the value of the standard intertwining operator on the spherical vector is the same than in the case of a non-Archimedean local field of characteristic $0$ with the same residue field \cite{Cs80, 3.1}. As the normalizing factors agree for both fields, the result is also true for $F_v$ of positive characteristic.

Suppose now that $v$ is an Archimedean place. Unfortunately, we could not find in the literature the assertion that in this case the action of the normalized intertwining operator on $\phi _{1,v}$ does not depend on $s$. So, we will give a proof of this hereafter. Assume for a moment already proved that the action of the normalized intertwining operator on $\phi_{1,v}$ does not depend on $s$. To compute the value in $(2)$, we can then assume $s=0$ and $\chi_v=1$, as $\chi _v$ is unramified by assumption and $r_{\alpha ,v}(s,\chi _{s_0\ti{\alpha }},\psi _{F_v})=r_{\alpha ,v}(s+s_0,1,\psi _{F_v})$, the $\epsilon $-factor being independent of $s$ in the Archimedean case. Then, we can apply \cite{Sh85, Theorem 4.1}, remarking that the trivial representation of the maximal torus induces irreducible for $SL_2$ and $SU(2,1)$ \cite{Ko69} (this is enough as the intertwining operator factors through the simply connected covering of the derived group \cite{Sh85, p. 978}) and using \cite{Sh85, Lemma 4.2} (and its proof) that shows $L(s,1,r_{\alpha }^{\vee })=L(s,1,r_{\alpha })$ (the assumptions of the lemma are satisfied, as after our scaling we have $\sigma =1$).

Let us now prove that the action of the normalized intertwining operator on $\phi _{1,v}$ does not depend on $s$ for $v$ an Archimedean place. Assume first that $v$ is real Archimedean. By the remark just above, one can restrict to $\ti{G^{der}}$, and it is enough to prove the assumption for each factor of $\ti{G^{der}}$ over $F_v$, so that we can assume $\ti{G^{der}}$ is isomorphic either to  $Res_{F'_{v'}/F_v}SL_2$ or to $Res_{F'_{v'}/F_v}SU(2,1)$ over $F_v$. As above, by applying a scaling, we can reduce to the case $\chi_v=1$. If $\ti{G^{der}}$ is of relative rank $1$ over $F_v$, then one has either $\ti{G^{der}}=SL_2$, $\ti{G^{der}}=Res_{\Bbb C/\Bbb R}SL_2$ or $\ti{G^{der}}\simeq SU(2,1)$.
By \cite{Kn86, prop. 7.6, p. 180} and \cite{Sch71, prop. 2.5} (remarking that $\langle s\ti{\alpha },\alpha ^{\vee }\rangle=4s$ in the relative root system for the quasi-split special unitary group \cite{Sh88} or \cite{Sh10, p.10, between 1.2.4 and 1.2.5}, and that $(p,q)=(1,0)$ if $\ti{G^{der}}\simeq SL_2$, $(p,q)=(2,0)$ if $\ti{G^{der}}\simeq Res_{\Bbb C/\Bbb R}SL_2$ and $(p,q)=(2,1)$ if $\ti{G^{der}}\simeq SU(2,1)$ in the notations of \cite{Sch71, Kn86}), one has, with a suitable normalization of the measures,
$$(A _v(s\ti{\alpha }, 1, \ti{w_{\alpha }})\phi_{1,v})(1)=\cases {\Gamma(1)\over\Gamma({1\over 2})}{\Gamma({s\over 2})\over\Gamma({s+1\over 2})}, &\hbox{\rm if $\ti{G^{der}}\simeq SL_2$;} \cr {\Gamma(2)\over\Gamma(1)}{\Gamma(s)\over\Gamma(s+1)}, & \hbox{\rm if $\ti{G^{der}}\simeq Res_{\Bbb C/\Bbb R}SL_2$;} \cr  {\Gamma(3)\over\Gamma({3\over 2})}{\Gamma(2s)\Gamma(s+{1\over 2})\over\Gamma(2s+1)\Gamma(s+1)}, & \hbox{\rm if $\ti{G^{der}}\simeq SU(2,1)$,}\cr \endcases\eqno{(3)}$$
where $\Gamma $ is the usual $\Gamma $-function. In addition, changing the measure does change the value of the integral only by a constant that does not depend on $s$ \cite{Sch71, (2.2.8)}. As we only want to show that the action of the normalized intertwining operator on $\phi _{1,v}$ does not depend on $s$, the precise measure does not matter at this moment.

If $\ti{G^{der}}=SL_2$ or $\ti{G^{der}}=Res_{\Bbb C/\Bbb R}SL_2$, then this agrees well, up to a constant, with the reciprocal of the normalizing factor from the Langlands-Shahidi theory which is in this case given by the Artin local factors, remarking that $\langle s\ti{\alpha },\alpha^{\vee }\rangle =2s$ in the second case. By the product formula for the intertwining operator and the above result for the standard intertwining operator in the relative rank one case, this takes also care of the case $\ti{G^{der}}\simeq Res_{F'_{v'}/F_v}SU(2,1)$ if $SU(2,1)$ splits over $F'_{v'}$.

Let us now suppose that $\ti{G^{der}}=SU(2,1)$ and that $SU(2,1)$ does not split over $\Bbb R$. Then, $\eta_{E'/F',v}$ is the $sgn$-character. The local factors from the Langlands-Shahidi theory associated to $r_{\alpha,1}^{\vee }$ are then $L_{\Bbb C}(s)$ and $\epsilon_{\Bbb C}(s,\psi _{F_v})$ which is constant, where $L_{\Bbb C}$ and $\epsilon _{\Bbb C}$ are the standard complex Artin local factors. Those associated to $r_{\alpha ,2}^{\vee }$, give the Asai-$L$-function. As $\Bbb C/\Bbb R$ corresponds by local class field theory to the sign character, one gets $L_{\Bbb R}(s,sign )$, while $\epsilon_{\Bbb R}(s, sign,\psi _{F_v})$ is also constant.

In explicit terms, $L_{\Bbb C}(s)=2(2\pi )^{-s}\Gamma(s)$ and $L_{\Bbb R}(s,sgn)=\pi^{-{s+1\over 2}}\Gamma({s+1\over 2})$. It follows that the normalizing factor $${\epsilon _v(s,1,r^{\vee }_1,\psi _{F_v})L_v(1+s,1,r^{\vee }_1)\over L_v(s,1,r^{\vee }_1)}{\epsilon _v(2s,1,r^{\vee }_2,\psi _{F_v})L_v(1+2s,1,r^{\vee }_2)\over L_v(2s,1,r^{\vee }_2)}$$
is up to a constant equal to $${\Gamma(s+1)\Gamma(s+1)\over\Gamma(s)\Gamma(s+{1\over 2})}.$$

As, by the Legendre duplication formula, one has $\Gamma(2s)={2^{2s-1}\over\sqrt{\pi }}\Gamma(s)\Gamma(s+{1\over 2})$ and in particular $\Gamma(2s+1)={2^{2s}\over\sqrt{\pi }}\Gamma(s+{1\over 2})\Gamma (s+1)$, one sees that this equals, upto a constant, the reciprocal of ${\Gamma(2s)\Gamma(s+{1\over 2})\over\Gamma(2s+1)\Gamma(s+1)}$. By the above formula $(3)$ for the value of the standard intertwining operator on the spherical vector, it follows that the normalized intertwining operator acts also by a constant if $\ti{G^{der}}=SU(2,1)$ and $SU(2,1)$ does not split over $\Bbb R$..

Finally, suppose $v$ is a complex Archimedean place of $F$. Then, $G_v$ splits necessarily over $F_v$ and one is reduced to $\ti{G^{der}}(F_v)=SL_2(\Bbb C)$ or $\ti{G^{der}}(F_v)\simeq SL_3(\Bbb C)$. If $\ti{G^{der}}(F_v)=SL_2(\Bbb C)$, then the value of the standard intertwining operator on the spherical vector has already been given above and the local factors are Hecke ones, also described above. One sees from that that the value of the normalized intertwining operator on the spherical vector $\phi _{1,v}$ is a constant function. By the product formula for the normalized intertwining operator, the same follows if $\ti{G^{der}}(F_v)\simeq SL_3(\Bbb C)$.\hfill{\fin 2}

\null{\bf Proposition:} (with the notations and assumptions of the theorem) \it Let $\Sigma _i$ be an irreducible component of the relative root system $\Sigma $ and denote by $G_{\Sigma_i}$ the simple factor of $G$ corresponding to $\Sigma _i$. Then, there is an everywhere unramified extension $F'/F$ and an absolutely simple quasi-split group  $G_i$ over $F'$, such that $G_{\Sigma_i}$ is isomorphic to $Res_{F'/F}G_i$. Denote by $d'$ the degree of the extension $F'/F$.

If $G_i$ splits over $F'$, then one has $d_{\alpha }=d'$ for all $\alpha\in\Sigma _i$;

If $G_i$ does not split over $F'$, then the simply connected covering $\ti{G^{der}_i}$ of $G_i$ is one of the following: $SU(n,n+1)$, $SU(n,n)$, $Spin_{2n}^-$, $^3D_4$ or $^2E_6$. One has

- for $\ti{G^{der}_i}=SU(n,n+1)$ (resp. $SU(n,n)$), $d_{\alpha }=2d'$ if $\alpha $ is a long (resp. short) root and $d_{\alpha }=d'$ if $\alpha $ is a short (resp. long) root;

- for $\ti{G^{der}_i}=Spin_{2n}^-$, $d_{\alpha }=2d'$ if $\alpha $ is a short root and $d_{\alpha }=d'$ if $\alpha $ is a long root;

- for $\ti{G^{der}_i}=\ ^3D_4$, $d_{\alpha }=3d'$ if $\alpha $ is a long root and $d_{\alpha }=d'$ if $\alpha $ is a short root;

- for $\ti{G^{der}_i}=\ ^2E_6$, $d_{\alpha }=2d'$ if $\alpha $ is a short root and $d_{\alpha }=d'$ if $\alpha $ is a long root.

\null Proof: \rm
This is obvious in the split case. In the non-split case, this can be read off for example from \cite{KK11} if $G_i$ is simply connected. In general, it is enough to remark that, for $\alpha\in\Sigma _i$, the simply connected covering of $(\ti{G^{der}_i})_{\alpha }$ is isomorphic to $\ti{G^{der}_{i,\alpha }}$ (where here the index $\alpha $ denotes the minimal Levis subgroup of semisimple rank $1$ associated to $\alpha $ in respectively $\ti{G^{der}_i}$ and $G_i$ to distinguish between both). This follows from the fact that the simply connected covering map $\ti{G^{der}_i}\rightarrow G_i^{der}$ induces an $F$-isogeny $(\ti{G^{der}_i})_{\alpha }^{der}\rightarrow G_{i,\alpha }^{der}$. As both groups are semi-simple of relative rank one, their simply connected coverings must be isomorphic over $F$. \hfill{\fin 2}

\null The following corollary will be useful when generalizing the results in \cite{DHO15}:

\null{\bf Corollary:} (with the notations and assumptions of the theorem) \it Let $\Sigma _p$ be an irreducible sub-root system of $\Sigma _{red}$ such that, for each $\alpha\in\Sigma _p$, the meromorphic function $r_{\alpha }(s,\chi\circ\underline{\alpha }^{\vee })$ has a pole $p_{\alpha }>0$. Then the set $\{(\alpha, p_{\alpha })\vert \alpha\in\Sigma _p\}$ is a subset of a set $\{(\alpha, p_{\alpha })\vert \alpha\in\Sigma _i\}$, where $\Sigma _i$ is an irreducible component of $\Sigma $, and either the $p_{\alpha }$ are all equal or $\Sigma _p$ is not a simply laced root system. In the latter case, the set of the $p_{\alpha }$ contains two elements with ratio given by

- for $\Sigma _i$ of type $B_n$ or $F_4$: $p_{\alpha }/p_{\beta }=2$ if $\alpha $ is a short and $\beta $ a long root;

- for $\Sigma _i$ of type $C_n$: $p_{\alpha }/p_{\beta }=2$ if $\alpha $ is a long and $\beta $ a short root;

- for $\Sigma _i$ of type $G_2$: $p_{\alpha }/p_{\beta }=3$ if $\alpha $ is a long and $\beta $ a short root.\rm

\null
\null
\centerline{\sc Acknowledgement}

\null
The material in this note was part of a joint project with
Marcelo De Martino and Eric Opdam. I thank both of them for allowing me to
publish this note separately. I thank Eric Opdam for useful discussions,
and for correcting an error in an earlier version. A special thanks to Prof.
F. Shahidi for an intense exchange on his work which lead to two footnotes and a remark.

\null
\null

\enddocument